\documentclass{article}
\usepackage{amsmath}
\usepackage{amssymb}

\renewcommand\a{\alpha}

\newcommand\C{{\mathbb C}}
\renewcommand\d{\delta}

\newcommand\g{\gamma}

\newcommand\op[1]{\mathop{\rm #1}\nolimits}
\newcommand\p{\partial}
\newcommand\po{$\hspace{-5pt}{\bf .}$ }
\newcommand\qed{\phantom{\underline{y}}\hfill\hfill$\square$}
\newcommand\R{{\mathbb R}}
\renewcommand\t{\times}

\newcommand\ti{\tilde}
\newcommand\ve{\varepsilon}
\newcommand\vp{\varphi}

\newcommand\we{\wedge}

\newcommand\1{{\bf 1}}

\newcommand\bib[1]{\bibitem[#1]{#1}}
\newtheorem{theorem}{Theorem}
\newtheorem{theo}{Theorem}

\newtheorem{rk}[theo]{Remark.}
\newenvironment{proof}[1][Proof]{\textbf{#1.} }{\qed}

\newenvironment{ex}{\trivlist \item[\hskip \labelsep{\bf Example.}]}{\endtrivlist}

\begin{document}

\title{{Deformation of big
pseudoholomorphic disks and application to the Hanh pseudonorm}}
\author{\textsc{Boris Kruglikov}
\smallskip \\
\textsf{Mat-Stat. Dept., University of Troms\o, Norway}
\smallskip \\
 \small kruglikov@math.uit.no}
 \date{}
\maketitle

\vspace{-20pt}
 \begin{abstract}
We simplify proof of the theorem that close to any
pseudoholomorphic disk there passes a pseudoholomorphic disk of
arbitrary close size with any pre-described sufficiently close
direction. We apply these results to the Kobayashi and Hanh
pseudodistances. It is shown they coincide in dimensions higher
than four. The result is new even in the complex case.
 \end{abstract}


We aim here to prove the following statement, which was proved by
another (analogous to the approach of \cite{NW}) and more
complicated method in \cite{K1}.

 \begin{theorem}\po
Let $(M^{2n},J)$ be an almost complex manifold and
 $$
f_0:(D_R,i)\to(M,J),\qquad(f_0)_*(0)e=v_0\ne0,
 $$
be a pseudoholomorphic disk. Here $e=1$ is the unit vector at
$0\in\C$. For every $\ve>0$ there exists a neighborhood ${\cal
V}_\ve (v_0)$ of the vector $v_0\in TM$ such that for each
$v\in\cal{V}_\ve$ there is a bit smaller pseudoholomorphic disk
 $$
f:(D_{R-\ve},i)\to(M,J),\qquad f_*(0)e=v.
 $$
The approximating curve $f$ can be embedded/immersed if such is
the curve $f_0$.
 \end{theorem}

This theorem was used in \cite{K1} for the proof of equivalence of
two definitions of Kobayashi pseudodistance $d_M$ in almost
complex category. In the second $d_M$ is associated via path
integration to the Kobayashi-Royden pseudonorm:
 $$
F_M(v)=\inf\{1/r\ |\ f:(D_r,i)\to (M,J),\ f_*(0)e=v\},\ v\in TM.
 $$
The above theorem assures $F_M$ to be upper semicontinous,
implying that
 $$
d_M(x,y)=\inf\{\int\g^*F_M\,|\,\g:[0,1]\to M, \g(0)=x,\g(1)=y\}.
 $$
is well-defined.

Moreover, since an embedded disk can always be perturbed to
embedded we prove simultaneously the main properties of the Hanh
pseudonorm $S_M(v)$, which is defined by the same formula as $F_M$
with an additional requirement on $f$ to be injective. This
pseudometric generates a pseudodistance via path integration, like
$F_M$ generates $d_M$, and this coincides (cf. \cite{K1}) with the
distance
 $$
h_M(x,y)=\op{inf}\sum_{k=1}^{m}d(z_k,w_k),
 $$
defined via injective chains $f_k:D_1\to(M^{2n},J)$,
$k=1,\dots,m$, $f_1(z_1)=p$, $f_m(w_m)=q$ and $f_k(w_k) =
f_{k+1}(z_{k+1})$, where $d$ is the Poincare metric on $D_1$.

Our approach to the close PH-curve existence result is similar to
that of \cite{S}, where the linearization of the structure $J$ was
made at a point. We linearize the structure along the disk and use
the reduction of the almost complex problem to a complex one via
the Green operator:
 $$
T_r:C^k(D_r,\C^n)\to C^{k+1}(D_r,\C^n),\ \ g(z)\mapsto\frac1{2\pi
i}\iint_{D_r}\frac{g(z)}{\zeta-z}\, d\zeta\we d\bar\zeta.
 $$
It is continuous in the Sobolev and H\"older norms (\cite{V}) and
obeys the identities: $\bar\p T_r=\op{Id}$,\quad
$T_r\bar\p|_{C^{k+1}_0}=\op{Id}.$

 \begin{proof}
We study at first the case, when the curve is embedded. Let $U$ be
a neighborhood of the shrunk PH-curve $f_0(D_{R-\ve})$. We can
assume \cite{K1} the disk is standard
$f_0(D_{R-\ve})=D_{R-\ve}\times\{0\}^{n-1}\subset\C^n$ and the
almost complex structure $J:U\to\op{End}_\R(\C^n)$, $J^2=-\1$,
along it is the standard complex structure $J(z)=J_0$ for all
$z\in D_{R-\ve}$. The equation for $f$ to be pseudoholomorphic
reads:
 $$
\bar\p f+q_J(f)\p f=0,\qquad
q_J(z)=[J_0+J(z)]^{-1}\cdot[J_0-J(z)],
 $$
that due to the above properties is equivalent to
 $$
\bar\p h=0,\qquad h=[\op{Id}+T_{R-\ve}\circ q_J(f)\circ\p]\, (f).
 $$
For $k\in\R\setminus\mathbb{Z}$, $k>1$, consider the map
 $$
 \begin{array}{rrl}
\Phi: & \mathcal{J}\t C^{k+1}(D_{R-\ve};U) &
\hspace{-4pt}\longrightarrow\
C^{k+1}(D_{R-\ve};\C^n), \\
& (J,s) & \hspace{-4pt}\longmapsto\ [\op{Id}+T_{R-\ve}\circ
q_J(f_0+s)\circ\p]\, (f_0+s),
 \end{array}
 $$
where $\mathcal{J}$ is a neighborhood of the given almost complex
structure $J$ in $C^k$-topology. We consider $U$ as the total
space of the "normal bundle", with the sections being denoted by
$s$, so that every map $f\in C^{k+1}(D_{R-\ve};U)$, that is
$C^1$-close to $(f_0)|_{D_{R-\ve}}$, has a unique representation
$f=f_0+s$.

The map $\Phi_J=\Phi(J,\cdot)$ is $C^k$-smooth and satisfies:
$\Phi_J(0)=f_0$, $\Phi'_J(0)=\op{Id}$.
It has the Taylor decomposition (with $\|\cdot\|$ being the
$C^{k+1}$-norm):
 $$
\Phi_J(s)=f_0+s+o(\|s\|).
 $$
Therefore $\op{Im}\Phi_J$ contains a small neighborhood of the
curve $f_0$.

Let $Z=(a,v)\in T\C^n$ and $h_Z(z)=a+v z$ be the holomorphic disk
in $U$, $z\in D_{R-\ve}$. It is close to $f_0$ whenever $Z$ is
close to $Z_0=(0,(1,0,\dots,0))\in T\C^n$. Define
 \begin{equation*}
f_Z=f_0+\Phi^{-1}_J(h_Z).
 \end{equation*}
It is a $J$-holomorphic $(R-\ve)$-disk, which satisfies:
$f_Z-h_Z=o(|Z-Z_0|)$.

Consider the $C^k$-map $\Psi:\C^{2n}\to\C^{2n}$,
$Z\mapsto(f_Z(0),(f_Z)_*(0)e)$. Since the above estimate implies
$\Psi'(Z_0)=\op{Id}$, the map $\Psi(Z)$ is a local
$C^k$-diffeomor\-phism of a neighborhood of $Z_0$. In particular,
for every $Z=(a,v)$ sufficiently close to $Z_0$ there exists a
pair $\ti Z=(\a,\zeta)$ such that $\Psi(\ti Z)=Z$.

Now the obtained map $f=f_{\ti Z}\vphantom{\dfrac22}$ is
$C^1$-close to $f_0$ and so is embedded. It is also smooth due to
the usual elliptic regularity (\cite{NW,MD,S}). If $f_0$ is
immersed, the reasoning is the same for the neighborhood $U$
obtained via $f_0$ by the pull-back.

In the general case for the map $f_0:(D_{R-\ve},i)\to (M,J)$ we
consider the graph ${\hat f}_0:(D_{R-\ve},i)\to (D_{R-\ve}\t
M,\hat J=i\t J)$, which is injective and apply the part of the
statement already proved.
 \end{proof}

 \begin{rk}\hspace{-8pt}
The proof implies persistence of big pseudoholomorphic disks (with
an insignificant loss of size) under perturbation not only of the
initial point and direction, but also of the almost complex
structure $J$ (note the role of $\mathcal{J}$ above). This
generalizes theorems 1.7 and 3.1.1(ii) in \cite{MD} and \cite{S}
respectively.
 \end{rk}

The properties of the Kobayashi-Royden pseudometric for almost
complex manifolds was discussed in \cite{K1}. Let us consider the
non-integrable version $S_M$ of the Hanh pseudometric. It is known
(\cite{O}) that it coincides with the Kobayashi-Royden
pseudometric $F_M$ for domains $M\subset\C^n$ of dimension $n>2$.
We generalize this to the non-integrable case.

 \begin{theorem}\po
$S_M=F_M$ for almost complex manifolds $(M^{2n},J)$, $n>2$.
 \end{theorem}

 \begin{proof}
Since $S_M\ge F_M$, it is enough to show that whatever small
$\ve>0$ is, any pseudoholomorphic disk of radius $R>0$ can be
approximated by an injective pseudoholomorphic disk of radius
$R-\ve$ with the same initial direction.

We give at first a new simple proof of the mentioned above theorem
from \cite{O}. Let $M\subset\C^n$ be a domain and $f:D_R\to M$ be
a holomorphic map. Denote $f_W(z)=f(z)-w_2z^2-w_3z^3$, $z\in
D_{R-\ve}$, $W=(w_2,w_3)\in\C^{2n}$. For small $W$ the map has
still the image in $M$. Also note that $f_W(0)=f(0)$ and
$f_W'(0)=f'(0)$.

By the Sard's theorem a generic $w_1\in\C^n$ is outside the set
 $$
\Bigl\{\dfrac{f(z_1)-f(z_2)}{z_1^2-z_2^2}\,|\,z_1,z_2\in
D_{R-\ve}\Bigr\}\cup\Bigl\{\dfrac{f'(z)}{2z}\,|\,z\in
D_{R-\ve}\Bigr\}.
 $$
For such a choice the map $f_{w_2,0}$ is injective outside the
anti-diagonal $\{z_2=-z_1\}$. Note that regularity of the origin
is preserved. So, switching on $w_3$ being generic, we get the map
$f_{w_2,w_3}$ to be injective everywhere.

In other words, the Sard's theorem implies that the set of
$W=(w_2,w_3)$ for which $f_W$ is not injective has measure zero
and so a generic pair of small vectors $w_2,w_3\in\C^n$ defines
the required approximating disk $f_W(z)$.

In the general complex case we should shift along some holomorphic
vector fields. This is achieved by the graph-lift
construction and a Royden's lemma \cite{R} that an embedded
holomorphic disk, shrunk a bit, has a Stein neighborhood.

It is easier, however, to consider the general case of almost
complex manifolds $(M,J)$ and to deduce the statement for
integrable $J$ as a corollary.

Denote by $\pi:D_{R-\ve}\times M\to M$ the projection. As in
theorem 1, the graph-lift $\hat f_0:D_{R-\ve}\to D_{R-\ve}\times
M$ can be deformed to the family $\hat f_{\hat W}=\hat
f_0+\hat\Phi^{-1}_{\hat J}(g_{\hat W})$, where $g_{\hat W}=\hat
w_0+\hat w_1z-\hat w_2z^2-\hat w_3z^3$, $\hat W=(\hat w_2,\hat
w_3)$, $\hat w_j\in\C^{n+1}$ and $\hat w_0=\hat\vp_0(\hat w_2,\hat
w_3)$, $\hat w_1=\hat\vp_1(\hat w_2,\hat w_3)$ are some
$C^k$-smooth functions, close to zero, such that $(\hat f^\d_{\hat
W}(0),(\hat f^\d_{\hat W})_*(0)e)
=\hat Z_0=(0,(1,0,\dots,0))\in T\C^{n+1}$. We identify above $\hat
f_0$ with $g_{\hat0}$, the first coordinate disk, and its
neighborhood with a ball $B\subset\C^{n+1}$, equipped with the
structure $\hat J=i\times J$.

Similarly to the first proof we get: $\hat f_{\hat W}=g_{\hat
W}+\rho_{\hat W}$, where $\rho_{\hat W}=o(|\hat W|)$. Now $\hat
f_{\hat W}$ is an embedding if $\pi\hat f_{\hat W}(z_1)\ne\pi\hat
f_{\hat W}(z_2)$ for $z_1\ne z_2$ and $\p\pi\hat f_{\hat
W}(z)\ne0$. We consider only the first, more complicated,
injectivity condition. It's negation is equivalent to $g_{\hat
W}(z_1)-g_{\hat W}(z_2)=[\rho_{\hat W}]|_{z_1}^{z_2}+\zeta$,
$\zeta\in D$, or after the projection:
 $$
w_2(z_1+z_2)+w_3(z_1^2+z_1z_2+z_2^2)=w_1+\dfrac{\pi\rho_{\hat
W}(z_2)-\pi\rho_{\hat W}(z_1)}{z_2-z_1}
 $$
The last equation is never satisfied for a.e.\ small $W=(w_2,w_3)$
in $\C^{2n}$. Actually for $w_1=\pi\vp_1(\hat w_2,\hat w_3)$ the
r.h.s.\ is $o(|W|)$. Thus the claim follows from the Sard's
theorem, if at least one of the coefficients of $w_2$ and $w_3$ is
not small. Since
 \begin{multline*}
\hspace{-9pt}D\times D=
\\ \ \ =\bigl[U_{5\d}(z_1=z_2=0)\bigr]\cup\bigl[D\times D\setminus
U_\d(z_1=-z_2)\bigr]\cup\bigl[D\times D\setminus
U_\d(z_1=(-\tfrac12\pm\tfrac{\sqrt{3}}2)z_2)\bigr].
 \end{multline*}
and the regularity at $(0,0)$ is preserved under small
perturbation we may achieve injectivity away from the
anti-diagonal by the quadratic perturbation and the in its
neighborhood by a cubic one. This finishes the proof.
 \end{proof}
 \vspace{3pt}

For $n=1$, when almost complex structures are automatically
integrable, the equality $S_M=F_M$ is equivalent to
contractibility ($M$ being a disk). In the case of $\C$-dimension
$n=2$ the equality may fail to hold (however arguments of
theorem 2 show that $F_M$ coincides with the pseudonorm $\tilde
S_M$ obtained via immersed disks).

 \begin{ex}
Consider the map $\vp_\a:D_1\to\C^2$, $z\mapsto(z(\a z-1)^2,\a
z^2(\a z-1))$, $|\a|>1$. It has a unique self-intersection point
$\vp_\a(0)=\vp_\a(1/\a)=0$, which is transversal:
$\vp_\a'(0)=(1,0)$, $\vp_\a(1/\a)=(0,1)$, and so non-removable.
For a neighborhood $U$ of the image $\op{Im}(\vp_\a)$ the
pseudonorms $F_U$ and $S_U$ are different.
 \end{ex}

In the product case the pseudonorms $S_M$ and $F_M$ were compared
in \cite{J}. It is however unclear if we can majorize $S_M\le
c\cdot F_M$, with a constant depending on $M$, or more generally,
if Kobayashi and Hanh hyperbolicities ($d_M$ or resp.\ $S_M$ being
a metric) are equivalent. Of course, the former implies the
latter.

It was shown in \cite{KO} that contractible tame almost complex
domains are hyperbolic. In other cases the hyperbolicity may be
lost.

 \begin{ex}
Consider the Reeb foliation of $\R^3$ with the standard $T^2$ as a
leaf. This foliation propagates via parallel transports to
$\R^{2n}=\R^3\times\R^{2n-3}$, $n\ge2$, and there exists an almost
complex structure on $\R^{2n}$ such that the new foliation is
pseudoholomorphic. Every domain containing such a torus is not
tame and is not hyperbolic. Note that for $n=2$ only a curve of
genus one can be realized as a pseudoholomorphic submanifold in an
almost complex $(\R^{2n},J)$ (\cite{Mo}). For $n>2$ the sphere
$S^2$ can be realized as a pseudoholomorphic submanifold,
providing similarly a non-tame and non-hyperbolic domain in
$(\R^{2n},J)$ (\cite{KO}).
 \end{ex}

 \begin{rk}\hspace{-8pt}
The result of the last theorem shows that the analogy between
geodesics and pseudoholomorphic disks, instructive in many
respects (\cite{Mo}), is however limited. Though the Nijenhuis
tensor naturally plays the role of the curvature \cite{K2}, there
are no analogs for the conjugate points in complex time curves
theory.
 \end{rk}
%
%
%

\end{document}